\documentclass[12pt]{article}
\usepackage{graphicx}
\usepackage{amssymb,amsmath,amsfonts,amsthm,amscd,mathrsfs,amsthm}
\usepackage{epstopdf}
\usepackage{amsbsy}
\usepackage{psfrag}

\input xy
\xyoption{all}
\CompileMatrices
\newtheorem{theorem}[equation]{Theorem}
\newtheorem{lemma}[equation]{Lemma}

\newtheorem{corollary}[equation]{Corollary}

\theoremstyle{definition}
\newtheorem{definition}[equation]{Definition}

\theoremstyle{remark}
\newtheorem{remark}[equation]{Remark}
\newtheorem{remark/definition}[equation]{Remark-Definition}

\numberwithin{equation}{section}

\newcommand{\abs}[1]{| #1|}
\newcommand{\inv}[1]{\ensuremath{#1^{-1}}}
\newcommand{\norm}[2]{\ensuremath{\mathrm{N}}_{#1}(#2)}

\newcommand{\A}{\ensuremath{\mathbb{A}}}

\newcommand{\kssu}{\ensuremath{K(s;\sigma;u)}}
\newcommand{\N}{\ensuremath{\mathbb{N}}}
\newcommand{\Z}{\ensuremath{\mathbb{Z}}}

\newcommand{\ig}{\ensuremath{I[G]}}

\newcommand{\zg}{\ensuremath{\Z[G]}}

\newcommand{\gal}{\ensuremath{\mathrm{Gal}}}

\newcommand{\vbar}[1]{\ensuremath{v(#1)+\Gamma_F}}

\newcommand{\chara}{\ensuremath{\textrm{char}}}

\newcommand{\slashfrac}[2]{\ensuremath{\raise1ex\hbox{#1}\kern-.2em/\kern-.30em\lower1ex\hbox{#2}}}

\newcommand{\eqnref}[1]{(\ref{#1})}

\newcommand{\barm}{\ensuremath{\overline{m}}}
\newcommand{\barn}{\ensuremath{\overline{n}}}

\newcommand{\barw}{\ensuremath{\overline{w}}}

\newcommand{\sigm}{\ensuremath{\sigma^{\overline{m}}}}

\newcommand{\Br}{\ensuremath{\mathrm{Br}}}
\newcommand{\aut}{\ensuremath{\mathrm{Aut}}}
\newcommand{\oD}{\ensuremath{\overline{D}}}
\newcommand{\oF}{\ensuremath{\overline{F}}}
\newcommand{\oE}{\ensuremath{\overline{E}}}
\newcommand{\G}{\ensuremath{\langle \sigma_1\rangle \times \ldots \times \langle \sigma_r\rangle}}
\newcommand{\ou}{\ensuremath{\overline{u}}}
\newcommand{\inn}{I}
\newcommand{\sign}{\ensuremath{\sigma^{\barn}}}
\newcommand{\oa}{\ensuremath{\overline{a}}}
\newcommand{\ob}{\ensuremath{\overline{b}}}
\newcommand{\cA}{\ensuremath{\mathcal{A}}}

\title{Prime to $p$ extensions of the generic abelian crossed product}
\author{Kelly McKinnie}
\begin{document}
\maketitle
\abstract
In this paper we prove that the noncyclic generic abelian crossed product $p$-algebras constructed by Amitsur and Saltman in \cite{AS} remain noncyclic after tensoring by any prime to $p$ extension of their centers.  We also prove that an example due to Saltman of an indecomposable generic abelian crossed product with exponent $p$ and degree $p^2$ 
remains indecomposable after any prime to $p$ extension.
\setcounter{section}{-1}
\section{Introduction}
Let $p$ be a prime.  A \emph{$p$-algebra} is a finite dimensional division algebra which is central over a field of characteristic $p$ and has index a power of $p$.  The question of whether or not every $p$-algebra is cyclic was first asked by Albert in \cite{Albert}.  It was a natural question to ask as he had proven in \cite{Albert-old} that every $p$-algebra is similar in the Brauer group of its center to a cyclic algebra.  Amitsur and Saltman answered the question in \cite{AS} by constructing, for any fixed prime $p$, generic abelian crossed product $p$-algebras which contain no $p$-power central elements, that is, non-central elements whose $p$-th power is central.  This proved the existence of noncyclic $p$-algebras of all degrees $p^n$, $n \geq 2$.  Credit for this problem is also due to Isaac Gordon who constructed a noncyclic 2-algebra of degree 4 in 1940.

Let $D$ be a finite dimension division algebra with center $F$ and let $E/F$ be an extension of fields.  We call $D_E=D\otimes_FE$ a \emph{prime to $p$ extension} of $D$ if the degree of $E/F$ is prime to $p$.  The main results of this paper, Theorem \ref{graded4} and Corollary \ref{cor89}, prove that there do not exist prime to $p$ extensions of the Amitsur-Saltman noncyclic generic abelian crossed products which are cyclic.  This is in the opposite direction of the well known result that every division algebra of prime degree does become cyclic after a prime to $p$ extension.  Other results regarding prime to $p$ extensions of $p$-power index division algebras can be found in \cite{RS}.  

The outline of this paper is as follows.  In section \ref{sec1} we determine precise conditions under which abelian crossed product algebras have homogeneous $p$-power central elements.  This condition, which is placed on the matrix defining the algebra, is called ``strong degeneracy'' \eqnref{defn2}.  In section \ref{sec1} we also show that non-degenerate and not strongly degenerate matrices remain so in any prime to $p$ extension \eqnref{theorem8}. 

In section \ref{sec2} we generalize the notions of degeneracy and strong degeneracy to valued division algebras which are semi-ramified with separable residue fields \eqnref{def1}.  This technique is used to prove that if $D/F$ is a valued division algebra with center $F$ a field of characteristic $p$ and $D/F$ is semi-ramified with residue field $\oD/\oF$ a separable extension and $D$ is not strongly degenerate, then $D$ contains no $p$-power central elements.  In \eqnref{theorem19} we show that all of these properties hold after a prime to $p$ extension. 

In section \ref{sec3} the main results of the paper, Theorem \ref{graded4} and Corollary \ref{cor89}, are obtained by applying the results of sections \ref{sec1} and \ref{sec2} to generic abelian crossed product $p$-algebras defined by not strongly degenerate matrices.  These are the same $p$-algebras that were proven to be noncyclic in \cite{AS} and we prove they remain noncyclic after any prime to $p$ extension of the center by showing they contain no $p$-power central elements.  In the final part of section \ref{sec3} we prove that in a special circumstance these algebras are indecomposable.

The work in this paper constitutes part of my thesis at the University of Texas at Austin.  Part of this work was completed while I held a continuing education fellowship from the University of Texas at Austin and a VIGRE graduate research fellowship.  I would like to thank my doctoral thesis advisor, Dr. David Saltman, for suggesting the problem to me and for many helpful conversations.  I would also like to thank Adrian Wadsworth for reading an earlier draft of this paper and for many valuable suggestions.

\section{Strong degeneracy}\label{sec1}
\subsection{Definition and characterization of strong degeneracy}
Let $G =\langle \sigma_1 \rangle \times \ldots \times \langle \sigma_r \rangle$, a finite abelian group with $\abs{\sigma_i}=n_i$.  Let $A$ be a central simple $F$-algebra with maximal Galois subfield $K$ and $\gal(K/F)=G$.  Then $A/F$ is well known to be an abelian crossed product. Following \cite{AS} the crossed product structure of $A$ can be given in the following way.  For each $1 \leq i \leq r$, choose $z_i \in A$ which extends by inner automorphism the action of $\sigma_i$ on $K$.  Given any $g \in G$ there exist $m_i \in \N$, $0 \leq m_i < n_i$, so that $g=\sigma_1^{m_1}\ldots \sigma_r^{m_r}$.  Setting $z_g=z_1^{m_1} \ldots z_r^{m_r}$ we have inner automorphism by $z_g$ on $A$ extends the action of $g$ on $K$.  Define $c:G \times G \to K^*$ by $z_g z_{g'}=c(g,g')z_{gg'}$.  The multiplicative associativity of $A$ implies that $c$ is a 2-cocycle.  Furthermore, $A$ is isomorphic to the crossed product $(K/F,G,c)$ generated over $K$ by the $z_g$.  

Set $u_{ij} = z_iz_j\inv{z_i}\inv{z_j}\in K$ and $b_i = z_i^{n_i}\in K$.  One can compute the 2-cocycle $c(g,g')$ using the elements $\{b_i\}_{i=1}^r$, $\{u_{ij}\}_{i,j=1}^r$ and their Galois conjugates.  Therefore we use the notation
\[A \cong (K/F,z_{\sigma},u,b):=\bigoplus_{\begin{array}{c}0\leq m_i <n_1\\ 1 \leq i \leq r\end{array}}K\cdot z_1^{m_1}\cdot\ldots\cdot z_r^{m_r}.\]
In $(K/F,z_{\sigma},u,b)$ we use the notation $z_{\sigma}$ because the $z_i$ were chosen with respect the the basis $\sigma = \{ \sigma_1, \ldots, \sigma_r\}$ of $G$, $u=(u_{ij})\in M_r(K^*)$ and $b=\{b_i\}_{i=1}^r\in (K^*)^r$.  In this algebra multiplication is gotten from the relations $z_iz_j=u_{ij}z_jz_i$ and $z_i^{n_1}=b_i$.  From \cite{AS}, Lemma 1.2, the $u_{ij}$ and $b_i$ are not independent.  In particular, they satisfy
\begin{gather}
u_{ii}=1; \hspace{.2in} u_{ji}=\inv{u_{ij}} \label{equation1}\\
\sigma_k(b_i)=\norm{i}{u_{ki}}b_i\label{equation2}\\
\norm{ik}{u_{ik}}=1 \label{equation3}\\
\sigma_i(u_{jk})\sigma_j(u_{ki})\sigma_k(u_{ij})=u_{jk}u_{ki}u_{ij},\label{equation4}
\end{gather}
where $\mathrm{N}_i$ (resp. $\mathrm{N}_{ij}$) is the norm map of the field extension $F/F^{\langle \sigma_i \rangle}$ (resp. $F/F^{\langle \sigma_i,\sigma_j \rangle}$).  Furthermore, by \cite{AS}, Theorem 1.3, if $K/F$ is a $G$-Galois extension and a set of elements $\{u_{ij},b_i\}_{i,j=1}^r$ satisfy relations (\ref{equation1}), (\ref{equation2}), and (\ref{equation4}), then there exists an abelian crossed product $(K/F,z_{\sigma},u,b)$ defined by the matrix $u$ and the vector $b$.

The properties of an abelian crossed product are governed by properties of the matrix $u = (u_{ij})\in M_r(K^*)$.  Let $G$ be as above with $r \geq 2$, and $K/F$ a $G$-Galois extension of fields.  Given an abelian crossed product $\Delta = \left( K/F,z_{\sigma},u,b\right)$, for $\barm=(m_1,\ldots,m_r) \in \N^r$ set $\sigma^{\barm}=\sigma_1^{m_1}\ldots \sigma_r^{m_r} \in G$ and $z^{\barm}=z_1^{m_1} \ldots z_r^{m_r}$.  For $\barm,\barn \in \N^r$, set \[u_{\barm,\barn}=z^{\barm}z^{\barn}\inv{(z^{\barm})}\inv{(z^{\barn})}.\]
In \cite{AS}, pg.81, the matrix $u = (u_{ij})$ is defined to be \emph{degenerate} if there are $\sigma^{\barm},\sigma^{\barn}$ in $G$ and $a,b \in K^*$ such that $\langle\sigma^{\barm},\sigma^{\barn}\rangle$, the group generated by $\sigma^{\barm} \hbox{ and } \sigma^{\barn}$, is not cyclic and $u_{\barm,\barn}=\sigma^{\barm}(a)\inv{a}\sigma^{\barn}(b)\inv{b}$.  In \cite{AS}, Theorem 3.2, non-degenerate matrices are used to construct generic abelian crossed product $p$-algebras which contain no non-central $p$-power central elements.  We introduce here a more restrictive definition of degeneracy which will determine precisely when generic abelian crossed product $p$-algebras have a nontrivial $p$-power central element (Lemma \ref{lemma21}).  For any $i \in \{1, \ldots, r\}$ and any $\barm \in \N^r$, set $u_{i,\barm}=u_{\overline{e}_i,\barm}$, where $\overline{e}_i$ is the standard $i$-th basis vector.
\begin{definition}
Let $\Delta=(K/F,z_{\sigma},u,b)$ be an abelian crossed product with $G=\G$ noncyclic of rank $r$.  The matrix $u=(u_{ij})$ is {\bf strongly degenerate} if there exists an element $\sigma^{\barm} \in G$ with prime order $q$, and a set of elements $l,x_1,\ldots,x_r \in K^*$ such that for all $1 \leq i \leq r$,
\begin{equation}u_{i,\barm}=\frac{\sigma^{\barm}(x_i)}{x_i}\frac{l}{\sigma_i(l)}. \label{equation19}
\end{equation}
\label{defn2}\end{definition}
Strong degeneracy and degeneracy are notions that will only be used when the group $G$ is noncyclic.  In that case, if a matrix is strongly degenerate it is also degenerate.  For, if $u$ is strongly degenerate then the subgroup generated by $\sigm$ and $\sigma_i$ in $G$ will not be cyclic for at least one $1 \leq i \leq r$.  The following key lemma characterizes the property of strong degeneracy for abelian crossed products in terms of the existence of homogeneous prime power central elements.  It is independent of the characteristic of the center.  Let $\bar{0} = (0,\ldots,0) \in \N^r$ be the zero vector.  
\begin{lemma}  Let $G$ be a noncyclic finite abelian group with basis $\sigma=\{\sigma_1,\ldots,\sigma_r\}$, $K/F$ a $G$-Galois extension, and $\Delta=(K/F,z_{\sigma},u,b)$ an abelian crossed product.  The matrix $u$ is strongly degenerate  if and only if there exists a prime $q$ such that $\Delta$ has a $q$-power central element of the form $lz^{\barm}$ for $l \in K^*$ and $\barm \in \N^r-\{\bar{0}\}$.\label{lemma22}
\end{lemma}
\proof  Assume $\Delta$ has a $q$-power central element of the form $lz^{\barm}$ with $l \in K^*$ and $\barm \ne 0$.  Since $(lz^{\barm})^q$ is central, $(\sigma^{\barm})^q=e$, the identity element of $G$.  Let $K^{\barm}$ be the subfield of $K$ fixed by $\sigma^{\barm}$.  $K/K^{\barm}$ is a cyclic degree $q$ extension and we denote its norm by $\textrm{N}_{\barm}$.  Since $(lz^{\barm})^q$ commutes with each $z_i$, we see \begin{equation}
z_i(lz^{\barm})^qz_i^{-1} = (lz^{\barm})^q=\norm{\barm}{l}(z^{\barm})^q.\label{equation33}
\end{equation}  
On the other hand, we calculate,
\begin{eqnarray}
z_i(lz^{\barm})^qz_i^{-1} &=&\sigma_i(\norm{\barm}{l})(u_{i,\barm}z^{\barm})^q\nonumber \\
&=&\norm{\barm}{\sigma_i(l)u_{i,\barm}}(z^{\barm})^q.\label{equation34}
\end{eqnarray}
Combining lines (\ref{equation33}) and (\ref{equation34}), we see that $\norm{\barm}{\sigma_i(l)l^{-1}u_{i,\barm}}=1$.  Therefore, by Hilbert's Theorem 90, there exists $x_i \in K^*$ such that
\[ u_{i,\barm} = \frac{\sigma^{\barm}(x_i)}{x_i}\frac{l}{\sigma_i(l)}.\]
Since this is true for each $1 \leq i \leq r$, the matrix $u$ is strongly degenerate.  Conversely, assume $u$ is strongly degenerate.  Let $u_{i,\barm} = \sigma^{\barm}(x_i)x_i^{-1}l\sigma_i(l^{-1})$ with $\sigm$ of prime order $q$ in $G$.  We show that $lz^{\barm}$ is $q$-power central.  Since $\sigma^{\barm}$ has order $q$ in $G$, it is clear that $lz^{\barm}$ is not central and $(lz^{\barm})^q$ commutes with all elements of $K$.  To check that $(lz^{\barm})^q$ is central, it suffices to show that it commutes with each $z_i$.  For all $1 \leq i \leq r$ we compute,
\begin{eqnarray*}
z_i(lz^{\barm})^qz_i^{-1} &=&\norm{\barm}{\sigma_i(l)u_{i,\barm}}(z^{\barm})^q\\
&=&\norm{\barm}{\sigma^{\barm}(x_i)x_i^{-1}l}(z^{\barm})^q\\
&=&\norm{\barm}{l}(z^{\barm})^q\\
&=&(lz^{\barm})^q.
\end{eqnarray*}
\endproof
\begin{corollary} Let $G$ be a noncyclic finite abelian group with basis $\sigma=\{\sigma_i\}_{i=1}^r$ and let $(K/F,z_{\sigma},u,b) \cong (K/F,w_{\sigma},v,d)$ be an isomorphism between abelian crossed products.  Then $u$ is  strongly degenerate if and only if $v$ is strongly degenerate.  Moreover, $u$ is degenerate if and only if $v$ is degenerate.\label{corollary3}
\end{corollary}
\proof Assume $u$ is a strongly degenerate matrix.  By Lemma \ref{lemma22} there exists a $q$-power central element of the form $lz^{\barm}$ for some $\sigm \in G$ of order $q \ne1$.  Let 
\[
\xymatrix{
\phi:(K/F,z_{\sigma},u,b) \ar[r]^{\sim}& (K/F,w_{\sigma},v,d)
}\]
be the isomorphism from the hypothesis.  We can assume $\phi$ is the identity on the maximal subfield $K$.  Set $a_{\barm}=\phi(z^{\barm})\inv{(w^{\barm})}\in K$.  Also, for each $1 \leq i \leq r$, set $a_{i} = \phi(z_i)w_i^{-1}\in K$.  It is easy to compute that $la_{\barm}w^{\barm}$ is a $q$-power central element in $(K/F,w_{\sigma},v,d)$ and therefore, by Lemma \ref{lemma22}, $v$ is strongly degenerate.  On the other hand assume $u$ is degenerate.  Then there exist elements $a,b \in K^*$ and elements $\sigm,\sigma^{\barn} \in G$ such that $\langle \sigm,\sigma^{\barn}\rangle$ is not cyclic and $u_{\barm,\barn}=\sigm(a)\inv{a}\sigma^{\barn}(b)\inv{b}$.  Set $a_{\barm}=\phi(z^{\barm})\inv{(w^{\barm})}\in K$ and $a_{\barn}=\phi(z^{\barn})\inv{(w^{\barn})}\in K$.  Then it is easy to check that $ba_{\barm}w^{\barm}$ and $\inv{a}a_{\barn}w^{\barn}$ commute and that this is equivalent to $v$ being degenerate.
\endproof

\subsection{Strong degeneracy in prime to $p$ extensions}
The main result in this section, Theorem \ref{theorem8}, shows that if $u$ is not strongly degenerate in the extension $K/F$ with $[K:F]=p^s$, some $s \geq 1$, then $u$ remains not strongly degenerate in the extension $KE/E$ for $E$ any prime to $p$ extension of $F$.  We first prove a preliminary lemma and corollary.

\begin{lemma}
Let $\Delta=(K/F,z_{\sigma},u,b)=(K/F,G,c)$ be any abelian crossed product with $G$ noncyclic and $c:G\times G\to K^*$ a 2-cocycle.  Then $(K/F,G,c^t)=(K/F,w_{\sigma},u^t,b^t)$ where $u^t=(u_{ij}^t)$ and $b^t=\{b_i^t\}$ for any integer $t >0$. \label{lemma9}
\end{lemma}
\proof
Let $B=(K/F,G,c^t)$ be generated over $K$ by generators $\{w_g\}_{g \in G}$ and relations $w_gw_h=(c(g,h))^t w_{gh}$ for all $g, h \in G$.  Set $w_{\sigma_i}=w_i$.  Then since $c$ is defined by $z_{\sigm}z_{\sign}=c(\sigm,\sign)z_{\sigm\sign}$ we have $w_{\sigm}=w^{\barm}$ and $B=(K/F,w_{\sigma},v,d)$ where
\begin{equation*}
\begin{split}
v_{ij}&=w_iw_j\inv{w_i}\inv{w_j}=w_iw_j\inv{w_{\sigma_i\sigma_j}}\inv{[w_jw_i\inv{w_{\sigma_j\sigma_i}}]} \\
 &= [c(\sigma_i,\sigma_j)c(\sigma_j,\sigma_i)^{-1}]^t=[z_iz_j\inv{z_i}\inv{z_j}]^t=u_{ij}^t.
\end{split}
\end{equation*}
Therefore the matrix $v=(v_{ij})$ is given by $v_{ij}=u_{ij}^t$.  We will denote this matrix as $u^t$, (which should not be confused with the $t$-th power of the matrix $u$).  Let $n_i=\abs{\sigma_i}$.  Computing $w_{i}^{n_i}$  we get,
\[
w_i^{n_i} = \prod_{j=1}^{n_i-1}c^t(\sigma_i^j,\sigma_i)w_{\sigma_i^{n_i}}
= \left(\prod_{j=1}^{n_i-1}c(\sigma_i^j,\sigma_i)\right)^t
=(z_i^{n_i})^t = b_i^t.
\]
Therefore, $d=b^t=\{b_i^t\}_{i=1}^r$ and $B = (K/F,w_{\sigma},u^t,b^t)$.
\endproof

\begin{corollary}
Let $\Delta=(K/F,z_{\sigma},u,b)$ be an abelian crossed product.  Assume $\Delta \sim \Delta^n$ in $\Br(F)$for some $n \in \Z$.  Then the matrix $u=(u_{ij})$ is strongly degenerate in $K/F$ if and only if $u^n=(u_{ij}^n)$ is strongly degenerate in $K/F$. Moreover, under the same hypothesis the matrix $u$ is degenerate in $K/F$ if and only if $u^n$ is degenerate in $K/F$.  \label{corollary1}
\end{corollary}
\proof  Let $\Delta=(K/F,G,c)$.  Then $\Delta^n \sim (K/F,G,c^n)$ and, using Lemma \ref{lemma9}, 
\[(K/F,G,c^n) = (K/F,w_{\sigma},u^n,b^n).\]  
Thus, $\Delta \sim \Delta^n$ implies that $(K/F,z_{\sigma},u,b) \cong (K/F,w_{\sigma},u^n,b^n)$ since the degrees of the algebras over $F$ are equal.  Therefore, by Corollary \ref{corollary3}, $u$ is strongly degenerate if and only if $u^n$ is strongly degenerate and $u$ is degenerate if and only if $u^n$ is degenerate.  \endproof

We can now prove the main result of this section.
\begin{theorem}  Let $K/F$ be a finite abelian $p$-extension of any characteristic with $\gal(K/F)=G=\langle \sigma_1 \rangle \times\ldots\times\langle \sigma_r \rangle$, $r \geq 2$.  Let $\Delta=(K/F,z_{\sigma},u,b)$ be an abelian crossed product with $u$ not strongly degenerate (resp. not degenerate).  Then $u$ remains not strongly degenerate (resp. not degenerate) after any prime to $p$ extension of $F$. \label{theorem8}
\end{theorem}
\proof  Let $E/F$ be a prime to $p$ extension and identify $G$ with the Galois group of $KE/E$.  Assume by way of contradiction that there exists an element $\sigma^{\barm} \in G$ with order $p$ and a set of elements $l,x_1,\ldots,x_r \in (KE)^*$ such that for all $1 \leq i \leq r$,
\begin{equation}u_{i,\barm}=\frac{\sigma^{\barm}(x_i)}{x_i}\frac{l}{\sigma_i(l)}. \label{eqn43}
\end{equation}
Taking the norm, $N=N_{KE/K}$, from $KE$ to $K$ of both sides of \eqnref{eqn43} gives
\begin{eqnarray*}
\norm{}{u_{i,\barm}}&=&u_{i,\barm}^t \\
&=&\textrm{N}\left( \frac{\sigma^{\barm}(x_i)}{x_i}\frac{l}{\sigma_i(l)}\right)=\frac{\sigma^{\barm}(\norm{}{x_i})}{\norm{}{x_i}}\frac{\norm{}{l}}{\sigma_i(\norm{}{l})},
\end{eqnarray*}
where $t=[KE:K]$.  Setting $x'_i=\norm{}{x_i}$ and $l'=\norm{}{l}$ we get
\[u_{i,\barm}^t=\frac{\sigma^{\barm}(x'_i)}{x'_i}\frac{l'}{\sigma_i(l')} \,\,\hbox{ with } x'_i,l' \in K^*.\]
That is, the matrix $u^t=(u_{ij}^t)$ is strongly degenerate in $K/F$ where $t$ is prime to $p$.  Set $e=\hbox{exp}(\Delta)$.  Since $(t,p)=1$, we also have $(t,e)=1$.  Thus, there exists $k,l \in \Z$ such that $tk+el=1$.  In particular we have 
\begin{equation}\Delta = \Delta^{tk+el} \sim \Delta^{tk}\otimes_F \Delta^{el} \sim \Delta^{tk}.\label{eqn0}\end{equation}
Notice that the matrix $u^{tk}$ is strongly degenerate since $u^t$ is strongly degenerate.  Therefore, by \eqnref{eqn0} and applying Corollary \ref{corollary1}, we see that $u$ is strongly degenerate in $K/F$, a contradiction.  Thus, $u$ remains not strongly degenerate after any prime to $p$ extension.  If $u$ is only assumed to be non-degenerate in $K/F$, replace the definition of strong degeneracy in the above argument with regular degeneracy.  This gives that $u^{tk}$ is degenerate.  Using Corollary \ref{corollary1} again we see that $u$ is degenerate.  This is a contradiction.  \endproof

\section{Strong degeneracy in valued division algebras}\label{sec2}
In this section we use the fundamental homomorphism from \cite{JW}, which we recall below, to translate the notion of strong degeneracy on an abelian crossed product to any semi-ramified division algebra $D/F$ with $\oD/\oF$ a separable extension.  This translation allows us to study when such an algebra over a field of characteristic $p$ has a nontrivial $p$-power central element (Corollary \ref{corollary10}).
\subsection{Another characterization of strong degeneracy}
Let $\mathcal{D}(F)$ be the collection of division algebras which are finite dimensional and central over the field $F$.    If $D\in \mathcal{D}(F)$ is endowed with a valuation $v$ then we will denote $D$ by $(D,v)$.  $\Gamma_D$ and $\Gamma_F$ will be the \emph{value groups} of $D$ and $F$ respectively.  We will also use the notation $V_D=\{d \in D|v(d)\geq 0\}$ for the \emph{valuation ring} of $D$ and $U_D=\{d \in D|v(d)=0\}$ for the \emph{group of units} of $D$.  Recall a valued division algebra $D \in \mathcal{D}(F)$ is said to be \textit{semi-ramified} if $\overline{D}$, the residue division ring of $D$, is a field and $[\overline{D}:\overline{F}]=\abs{\Gamma_D:\Gamma_F}=\sqrt{[D:F]}$.  

Fix $(D,v) \in \mathcal{D}(F)$ with $(D,v)$ a semi-ramified division algebra and $\overline{D}/\overline{F}$ is a separable extension.  Under these hypotheses, by \cite{JW}, Proposition 1.7, the group homomorphism  
\begin{equation} \theta_{D}:\Gamma_D/\Gamma_F \to \gal(\overline{D}/\overline{F})\label{graded6}\end{equation}
is an isomorphism and in particular $\oD/\oF$ is an abelian Galois extension.  Recall for any $\pi \in D^*$, $\theta_D(v(\pi)+\Gamma_F)\in \gal(\overline{D}/\overline{F})$ is the automorphism on $\overline{D}$ induced by conjugation by $\pi$ on $V_D$, see \cite{JW} (1.6).  Set $G=\gal(\overline{D}/\overline{F})$ and choose a basis $\sigma = \{\sigma_i\}_{i=1}^r$ of the finite abelian group $G$ so that $G \cong \langle \sigma_1\rangle \times \ldots \times \langle \sigma_r \rangle$.  For $i=1,\ldots,r$ choose $\pi_i \in D$ such that $\theta_D(v(\pi_i)+\Gamma_F)=\sigma_i$.  Set $u_{ij}=\pi_i\pi_j\inv{\pi_i}\inv{\pi_j} \in U_D$ and set $\overline{u}_{ij}$ to be the image of $u_{ij}$ in the residue field $\overline{D}$.  Denote by $\overline{u}$ the matrix $(\overline{u}_{ij})\in M_r(\oD)$.  For any $\barm,\barn \in \N^r$, set $u_{\barm,\barn}=\pi^{\barm}\pi^{\barn}\inv{(\pi^{\barm})}\inv{(\pi^{\barn})} \in U_D$ and let $\overline{u}_{\barm,\barn}$ denote its image in $\overline{D}$.

Let $GD$ be the associated graded division algebra of the valued division algebra $D$ (\cite{Wadsworth-graded}, section 4).  By definition, $GD=\oplus_{\gamma \in \Gamma_D}GD_{\gamma}$ where $GD_{\gamma}$ is the residue ring $GD_{\gamma}=W^{\gamma}/W^{>\gamma}$ with $W^{\gamma}=\{d \in D^* | v(d)\geq \gamma\} \cup \{0\}$ and $W^{>\gamma}=\{d \in D^*| v(d)>\gamma\}\cup\{0\}$.  Let $g:D^* \to GD$ be the natural homomorphism taking a nonzero element of $D$ to its associated homogeneous element of $GD$.  Since $D$ satisfies the three properties $[\oD:\oF]\,|\Gamma_D:\Gamma_F| = [D:F]$, $\oD$ is separable over $\oF$ and $\chara(\oF) \nmid |\textrm{ker}(\theta_D)|$, by \cite{Boulagouaz2}, Corollary 4.4, $Z(GD)=GF$.  Here $GF$ is the associated graded field of $F$ and in general $GF$ is only known to be contained in $Z(GD)$.  Moreover, inner automorphism by $g(\pi)$ on $D_0=\overline{D}$, as an element of $\aut_{\overline{F}}(\overline{D})$, is equal to $\theta_D(v(\pi)+\Gamma_F) \in \gal(\overline{D}/\overline{F})$ (see \cite{Wadsworth-graded} (4.9)).  The following theorem characterizes degeneracy and strong degeneracy of the matrix $\overline{u}$ in the extension $\oD/\oF$ in terms of the behavior of elements of $D$ and can be viewed as a generalization of Lemma \ref{lemma22}.

\begin{theorem}
Let $(D,v) \in \mathcal{D}(F)$ be a valued division algebra, semi-ramified with $\oD$ separable over $\oF$ and $\gal(\oD/\oF)=G=\G$ a noncyclic group.  Choose $\{\pi_i\}_{i=1}^r \subset D$ satisfying $\theta_D(v(\pi_i)+\Gamma_F)=\sigma_i$ in the isomorphism $\theta_D:\Gamma_D/\Gamma_F \to \gal(\oD/\oF)$.  Set $u_{ij}=\pi_i\pi_j\inv{\pi_i}\inv{\pi_j}$.  
\begin{itemize}
\item[(i)] The matrix $\ou=(\ou_{ij})$ is strongly degenerate in $\oD/\oF$ if and only if there exists a prime power central, homogeneous element in $GD_{\gamma}$ with $\gamma \in \Gamma_D - \Gamma_F$.
\item[(ii)] The matrix $\ou=(\ou_{ij})$ is degenerate in $\oD/\oF$ if and only if there exist homogeneous elements $g(\alpha) \in GD_{\gamma}$ and $g(\beta) \in GD_{\epsilon}$ such that $\langle \theta_D(\gamma+\Gamma_F),\theta_D(\epsilon+\Gamma_F)\rangle$ is a noncyclic subgroup of $G$ and $g(\alpha)$ and $g(\beta)$ commute in $GD$.
\end{itemize}
\label{graded1}
\end{theorem}
\proof  For ($i$), let $\alpha \in GD_{\gamma}$ be a homogeneous element of $GD$ such that $\gamma \notin \Gamma_F$ and $\alpha^q \in GF=Z(GD)$ for some prime $q$.  Since $\theta_D$ is an isomorphism we can express $\alpha=\alpha_0 g(f) g(\pi^{\barm})$ for some $\alpha_0 \in D_0$, $f \in F$ and $\barm \in \N^r$, such that $v(f\pi^{\barm})=\gamma$.  Since $\alpha^q$ is central if and only if $(\alpha_0g(\pi)^{\barm})^q$ is, we will assume without loss of generality that $g(f)=1$.  Since $\alpha^q \in GF$, $\inn_{\alpha^q}(z)=\alpha^q z \inv{(\alpha^q)}=z$ for all $z \in D_0$.  And, since $\inn_{\alpha}|_{D_0}=\sigm$, this implies that $(\sigma^{\barm})^q=e$, the identity element of $G$.  Let $\oD^{\barm}$ be the subfield of $\oD$ fixed by $\sigm$.  $\oD/\oD^{\barm}$ is a nontrivial cyclic extension of degree $q$ since $\gamma \notin \Gamma_F$ implies that $\barm \ne \overline{0}$.  Since $\alpha^q$ commutes with each $g(\pi_i)$ for $i=1,\ldots,r$ and $\inn_{g(\pi_i)}|_{D_0}=\sigma_i$, 
\[\norm{\barm}{\alpha_0}(g(\pi)^{\barm})^q=\alpha^q=g(\pi_i)\alpha^q\inv{g(\pi_i)}=\norm{\barm}{\sigma_i(\alpha_0)\overline{u}_{i,\barm}}(g(\pi)^{\barm})^q.\]
That is, in $D_0=\oD$, $\norm{\barm}{\sigma_i(\alpha_0)\inv{\alpha_0}\overline{u}_{i,\barm}}=1$.  Therefore, by Hilbert's Theorem 90, for each $i=1, \ldots, r$, there exists an $x_i \in \oD$ such that $\sigma_i(\alpha_0)\inv{\alpha_0}\overline{u}_{i,\barm}=\sigma^{\barm}(x_i)\inv{x_i}$.  This proves that the matrix $u$ is strongly degenerate.

Conversely, assume there exists a choice of elements $\pi_i \in D$ such that the corresponding matrix $\overline{u}$ is strongly degenerate.  In particular, there exists a prime $q$ and an element $\sigma^{\barm} \in G$ with order $q$, and elements $l,x_1,\ldots,x_r \in \oD$ such that $\overline{u}_{i,\barm}=\sigma^{\barm}(x_i)\inv{x_i}l\inv{\sigma_i(l)}$.  We show that the homogeneous element $l g(\pi)^{\barm}$ is $q$-power central in $GD_{v(\pi^{\barm})} \subset GD$.  This will complete the proof since $\sigm$ has order $q$ in $G$ and therefore $v(\pi^{\barm})\in \Gamma_D-\Gamma_F$.  It is enough to show that $l g(\pi)^{\barm}$ commutes with homogeneous elements.  Again since $\theta_D$ is an isomorphism, every homogeneous element can be expressed as $\alpha_0g(f)g(\pi^{\barn})$ for some $\barn \in \N^r$.  Since $\sigm$ has order $q$, and $lg(\pi^{\barm})\alpha_0=\sigm(\alpha_0)lg(\pi^{\barm})$, $\alpha_0$ commutes with $(lg(\pi)^{\barm})^q$.  Moreover,  $(lg(\pi^{\barm}))^q$ commutes with $g(f)$, since $f \in F$ and therefore it is only left to show that $(lg(\pi)^{\barm})^q$ commutes with $g(\pi^{\barn})$.  For this it suffices to show that $(lg(\pi)^{\barm})^q$ commutes with $g(\pi_i)$ for each $i =1,\ldots,r$.
\begin{equation}
\begin{split}
g(\pi_i)(lg(\pi)^{\barm})^q\inv{g(\pi_i)} 
&=g(\pi_i) \norm{\barm}{l}(g(\pi)^{\barm})^q \inv{g(\pi_i)}\\
&=\norm{\barm}{\sigma_i(l)\overline{u}_{i,\barm}}(g(\pi)^{\barm})^q\\
&=\norm{\barm}{l}(g(\pi)^{\barm})^q\\
&=(lg(\pi)^{\barm})^q
\end{split}\label{eqn90}
\end{equation}
The third equality in \eqnref{eqn90} follows from setting $\ou_{i\barm}=\sigma^{\barm}(x_i)\inv{x_i}l\inv{\sigma_i(l)}$.  This shows that $lg(\pi)^{\barm}\in GD_{v(\pi^{\barm})}$ is a prime power central homogeneous element of $GD$.  Since $\theta_D(v(\pi^{\barm})+\Gamma_F) \ne e$, the identity element of $G$, $v(\pi^{\barm}) \notin \Gamma_F$.  

To prove ($ii$), assume the matrix $\ou$ is degenerate.  Then, there exist $\sigm,\sign \in G$ and $\overline{a},\overline{b}\in \oD$ so that $\langle \sigm,\sign\rangle$ is noncyclic and
\[\overline{u}_{\barm,\barn}=\frac{\sigm(\overline{a})}{\overline{a}}\frac{\sign(\overline{b})}{\overline{b}}.\]
It is a direct computation to see that $\ob g(\pi^{\barm})$ and $\oa^{-1} g(\pi^{\barn})$ commute:
\begin{equation}
\begin{split}
\ob g(\pi^{\barm})\oa^{-1} g(\pi^{\barn})&=\ob \sigm(\oa^{-1})\ou_{\barm,\barn}\sign(\ob^{-1})\oa\,\oa^{-1} g(\pi^{\barn})\ob g(\pi^{\barm})\\
&= \oa^{-1} g(\pi^{\barn})\ob g(\pi^{\barm}).
\end{split}\nonumber
\end{equation}
Conversely, Assume there exist homogeneous elements $g(\alpha) \in GD_{\gamma}$ and $g(\beta) \in GD_{\epsilon}$ so that $\langle \theta_D(\gamma+\Gamma_F),\theta_D(\epsilon+\Gamma_F)\rangle$ is noncyclic and $g(\alpha)$ and $g(\beta)$ commute.  Set $\theta_D(\gamma+\Gamma_F)=\sigm$ and $\theta_D(\epsilon+\Gamma_F)=\sign$.  We can express $g(\alpha)=\oa g(\pi^{\barm})g(f_1)$ and $g(\beta)=\ob g(\pi^{\barn})g(f_2)$ for some $f_1,f_2 \in F$ and $\overline{a},\overline{b} \in D_0$.  From the relation $g(\alpha)g(\beta)=g(\beta)g(\alpha)$ we get
\begin{equation}
\begin{split}
\overline{1}&=g(\alpha)g(\beta)g(\alpha^{-1})g(\beta^{-1})\\
&=\oa g(\pi^{\barm})\ob g(\pi^{\barn})(\oa g(\pi^{\barm}))^{-1}(\ob g(\pi^{\barn}))^{-1}\\
&=\oa \sigm(\ob)\ou_{\barm,\barn}\sign(\oa^{-1})\ob^{-1}.
\end{split}\nonumber
\end{equation}
Therefore the matrix $\ou$ is degenerate.\endproof

\begin{remark/definition}  Theorem \ref{graded1}($i$) shows that we may define a semi-ramified valued division algebra with separable residue fields to be {\bf strongly degenerate} if its associated graded division algebra has prime power central homogeneous elements in $GD_{\gamma}$ for some $\gamma \notin \Gamma_F$.  In particular, the strong degeneracy of any matrix $\overline{u}$ associated to a choice of $\{\pi_i\}_{i=1}^r$ does not depend on the choice of the elements $\{\pi_i\}$.  Similarly, we may define a semi-ramified valued division algebra with separable residue fields to be {\bf degenerate} if it satisfies the equivalent properties in Theorem \ref{graded1}($ii$).\label{def1}  
\end{remark/definition}
\begin{remark}  Using the fundamental homomorphism $\theta_D$ to distinguish elements of a finite dimensional, central, semi-ramified algebra was also done in \cite{Boulagouaz}, Proposition 2.2.  In this proposition an explicit form with homogeneous elements of the central localization of a central semi-ramified graded division algebra over a graded field is given which is analogous to the canonical elements of a generic abelian crossed product.
\end{remark}
\subsection{Connection with the $I\otimes N$ decomposition}
Let $(D,v) \in \mathcal{D}(F)$ be a semi-ramified division algebra with $\oD/\oF$ a separable noncyclic Galois extension with $\gal(\oD/\oF)\cong\langle \sigma_1\rangle \times \ldots \times \langle \sigma_r \rangle$.  Let $\ou$ be the matrix defined in Theorem \ref{graded1} associated to a choice of elements $\pi_i \in D^*$ with $\theta_D(\pi_i)=\sigma_i$.  The matrix $\overline{u}$ defines an abelian crossed product as follows.   Let $n_i \in \N$ be the order of $\sigma_i$.  Since $\theta_D$ is an isomorphism and $\theta_D(n_i\vbar{\pi_i}) = \sigma_i^{n_i}=e$, the identity element of $G$, we have $v(\pi_i^{n_i}) \in \Gamma_F$.  Choose $\{f_i\}_{i=1}^r \subset F$ such that $v(f_i) = v(\pi_i^{n_i})$ and set $\pi_i^{n_i}\inv{f_i} = b_i\in U_D$.  Set $\overline{b} = \{\overline{b}_i\}_{i=1}^r$.
 \begin{lemma}
 $\overline{u}$ and $\overline{b}$ satisfy conditions (\ref{equation1}), (\ref{equation2}) and (\ref{equation4}) in the field extension $\oD/\oF$.  In particular, there exists an abelian crossed product \[(\oD/\oF, z_{\sigma},\overline{u},\overline{b})\] defined by $\overline{u}$ and $\overline{b}$. \label{lemma5}
 \end{lemma}
 \proof  Since $u_{ii} = 1$ and  $\inv{u_{ij}} = u_{ji}$, their images in $\oD$ satisfy condition \eqnref{equation1}.  To see that equation (\ref{equation2}) holds, note that by the definition of $\theta_D$ we have $\theta_D(\vbar{\pi_k}) = \sigma_k$.  Therefore $\sigma_k(\overline{b}_i)$ is equal to the image in $\oD$ of $\psi_{\pi_k}(b_i) = \pi_k b_i\inv{\pi_k}$.  Computing $ \pi_k b_i\inv{\pi_k}$ we get
 \begin{eqnarray}
 \pi_k b_i\inv{\pi_k} &=& \pi_k \pi_i^{n_i}\inv{f_i}\inv{\pi_k} =\pi_k \pi_i^{n_i}\inv{\pi_k}\inv{f_i} \nonumber \\
 &=&u_{ki}\pi_i\pi_k\pi_i^{n_i-1}\inv{\pi_k}\inv{f_i} = \ldots =u_{ki} \psi_{\pi_i}(u_{ki}) \ldots \psi_{\pi_i}^{n_i-1}(u_{ki})\pi_i^{n_i}\inv{f_i}. \nonumber
 \end{eqnarray}  
 Therefore, looking at the residues in $\oD$,
 $$\sigma_k(\overline{b}_i) = \overline{u}_{ki}\sigma_i(\overline{u}_{ki}) \ldots \sigma_i^{n_i-1}(\overline{u}_{ki})\overline{b}_i = \norm{i}{\overline{u}_{ki}}\overline{b}_i,$$
 proving that $\overline{u}$ and $\overline{b}$ satisfy (\ref{equation2}).  To prove that $\overline{u}$ satisfies \eqnref{equation4}, consider the following two different ways to compute the image of the element $\pi_k \pi_i \pi_j \inv{\pi_k} \inv{\pi_j} \inv{\pi_i}$ in $\oD$.
 \begin{eqnarray}
 \overline{\pi_k \pi_i \pi_j \inv{\pi_k} \inv{\pi_j} \inv{\pi_i}} &=& \overline{u_{ki} \pi_i u_{kj} \inv{\pi_i}} = \overline{u}_{ki} \sigma_i(\overline{u}_{kj})\label{equation6} \\
 \overline{\pi_k \pi_i \pi_j \inv{\pi_k} \inv{\pi_j} \inv{\pi_i}} &=& \overline{\pi_k u_{ij} \pi_j \pi_i \inv{\pi_k}\inv{\pi_j}\inv{\pi_i}} = \sigma_k(\overline{u}_{ij})\overline{u}_{kj}\overline{\pi_ju_{ki}\pi_i\inv{\pi_j}\inv{\pi_i}} \nonumber \\
 &=& \sigma_k(\overline{u}_{ij}) \overline{u}_{kj} \sigma_j(\overline{u}_{ki})\overline{u}_{ji}  \label{equation7}
 \end{eqnarray}
 Combining lines (\ref{equation6}) and (\ref{equation7}), we see that $\overline{u}$ satisfies \eqnref{equation4}.  Therefore, by \cite{AS} Theorem 1.3, there exists an abelian crossed product $(\oD/\oF,z_{\sigma},\overline{u},\overline{b})$. \endproof
\begin{remark}  The abelian crossed product in Lemma \ref{lemma5} defined by the matrix $\overline{u}$  is the same as one constructed in the proof of \cite{JW}, Theorem 5.6 in the case of an inertially split division algebra.  In that proof, for any inertially split division algebra, $D/F$, Jacob and Wadsworth construct a decomposition of $D$, so that $D$ is the underlying division algebra of $I \otimes N$, where $I$ is an inertial division algebra and $N$ is a nicely semi-ramified division algebra.  The abelian crossed product from Lemma \ref{lemma5} is precisely the residue division algebra of the inertial $I$ constructed in the decomposition $I \otimes N$ in the case when $D$ is both semi-ramified and inertially split.  Note that if we assumed our semi-ramified division algebras with separable residue fields were defined over a Henselian valued field it would be an inertially split division algebra.  This follows because there would exist an inertial lift of $\oD$ in $D$ which was a maximal subfield (see \cite{JW} page 148).  \end{remark}

\subsection{$p$-Power central elements in $p$-algebras}
In this section we state and prove the main corollary that we have been working up to, Corollary 2.7, regarding $p$-algebras with no $p$-power central elements.  We start with a lemma which was pointed out to me by Adrian Wadsworth.  
\begin{lemma} Let $(D,v) \in \mathcal{D}(F)$ be a valued division algebra with $F$ a field of characteristic $p$.  Assume
\begin{itemize}
\item[(i)]  $D$ is defectless over $F$,
\item[(ii)]  $\oD$ contains no proper purely inseparable extension of $\oF$, and
\item[(iii)]  if $a \in D$ with $a^p \in F$, then $v(a)\in \Gamma_F$, the value group of $F$.
\end{itemize}
Then for $a \in D$ with $a^p \in F$, $a \in F$.
\label{lemma33}
\end{lemma}
\proof  Assume by way of contradiction that there exists an $a \in D-F$ with $a^p \in F$.  Then $K=F(a)$ is a purely inseparable field extension of $F$ in $D$.  Since $D$ is defectless over $F$ by ($i$), $K$ is also defectless over $F$, that is, $[K:F] = [\overline{K}:\oF]\cdot [\Gamma_K:\Gamma_F]$, where $\Gamma_K$ is the value group of $K$.  By ($iii$), $D$ does not contain any totally ramified, purely inseparable field extensions, hence $\Gamma_K = \Gamma_F$.  Therefore, $\overline{K}/\oF$ is an extension of degree $p$.  However, the residue field of a purely inseparable extension is itself purely inseparable, hence $\overline{K}/\oF$ is a proper purely inseparable field extension of $\oF$, contained in $\oD$, a contradiction to ($ii$).  \endproof

We can now state our main result on the existence of $p$-power central elements in a $p$-algebra.

\begin{corollary}
Let $(D,v) \in \mathcal{D}(F)$ be a valued division algebra which is semi-ramified with $\oD/\oF$ a separable extension.  Assume $G=\gal(\oD/\oF)$ is a noncyclic $p$-group, $D$ is not strongly degenerate and $F$ is a field of characteristic $p$.  Then $D$ has no non-central $p$-power central elements.  
\label{corollary10}
\end{corollary}
\proof  By assumption, $D$ satisfies ($i$) and ($ii$) of Lemma \ref{lemma33}.  Let $x\in D$ be a $p$-power central element.  Then $g(x) \in GD$ is a homogeneous, $p$-power central element. Therefore, by Theorem \ref{graded1}($i$), $v(x) \in \Gamma_F$.  Therefore $D$ satisfies ($iii$) of Lemma \ref{lemma33} and we are done.  \endproof

\begin{remark}This theorem compares to \cite{Boulagouaz}, Lemma 3.1 and Theorem 3.2.  In \cite{Boulagouaz}, Boulagouaz and Mounirh use semi-ramified, central graded division algebras to get their result on prime power central elements.  It will be convenient for us to have the result over valued fields because we will apply \eqnref{corollary10} to extensions of Henselian valued fields.
\end{remark}

\subsection{Prime to $p$ extensions of division algebras}
Throughout this section let $(D,v) \in \mathcal{D}(F)$  with $\mathrm{deg}(D)=p^n$, $n\geq 1$ and no assumption on the characteristic of $F$.  Let $E/F$ be any prime to $p$ extension of $F=Z(D)$ and let $D_E = D\otimes_F E$ be the division algebra gotten by extending scalars to $E$
.  The next easy lemma shows that the properties of being a valued, semi-ramified division algebra with separable residue fields are preserved by prime to $p$ extensions.
\begin{lemma}Let $(D,v)\in \mathcal{D}(F)$ be a degree $p^n$ semi-ramified division algebra with $\oD/\oF$ a separable extension.  Let $E/F$ be a prime to $p$ extension of $F$.  Then $D_E=D\otimes_FE$ is a semi-ramified valued division algebra with $\overline{D_E}/\oE$ a separable extension.
\label{lemma34}
\end{lemma}

\proof Let $w$ be an extension of the valuation $v$ on $F$ to $E$.  Then $(D,v)$ and $(E,w)$ are valued division rings with $F=Z(D)\subseteq E$.  Moreover, $v|_F=w|_F$ by assumption and $[D:F]<\infty$.    The pair $(D,v)$ and $(E,w)$ also satisfy
\begin{enumerate} 
\item $D$ is defectless over $F$,
\item $\Gamma_D \cap \Gamma_E=\Gamma_F$, and
\item $\oD \otimes_{\oF}\oE$ is a division ring (in fact it is a field).
\end{enumerate}
Therefore, by \cite{Morandi}, Theorem 1, $D_E=D\otimes_FE$ is a division ring with the following properties.  There exists a valuation on $D_E$ extending both $v$ and $w$, $\Gamma_{D_E}=\Gamma_D+\Gamma_E$ and $\overline{D_E}\cong \oD \otimes_{\oF}\oE$.  In particular, $\overline{D_E}$ is the field join of $\oD$ and $\oE$, which is a separable extension of $\oE$ and $[\overline{D_E}:\oE]=[\oD:\oF]$.  As for the relative value groups, we compute that they satisfy,
\[\abs{\Gamma_{D_E}:\Gamma_E} = \frac{\abs{\Gamma_{D_E}:\Gamma_F}}{\abs{\Gamma_E:\Gamma_F}} = \frac{ \abs{\Gamma_{D_E}:\Gamma_D} \cdot \abs{\Gamma_D:\Gamma_F}}{\abs{\Gamma_E:\Gamma_F}}.\]
Since $|\Gamma_E:\Gamma_F|$ divides $[E:F]$, it is prime to $p$, and hence must divide $ \mid \Gamma_{D_E}:\Gamma_D\mid$.  Therefore, 
\begin{equation}\mid \Gamma_{D_E}:\Gamma_E\mid \geq \mid \Gamma_D:\Gamma_F\mid.\label{eqn48}\end{equation}
This inequality is an equality by the following argument.  Since $D_E$ is a valued division algebra it satisfies the ``fundamental inequality'' $[D_E:E]\geq [\overline{D_E}:\oE]\, |\Gamma_{D_E}:\Gamma_E|$.  Therefore, if \eqnref{eqn48} is strict inequality then
\[[D:F]=[D_E:E]\geq [\overline{D_E}:\oE]\abs{\Gamma_{D_E}:\Gamma_E}>[\oD:\oF]\abs{\Gamma_D:\Gamma_F},\]
which contradicts the defectlessness of $D$.  Therefore, $D_E$ is a valued, semi-ramified division algebra with separable residue fields.\endproof

As a consequence of Lemma \ref{lemma34} we can formulate and prove the following theorem.
\begin{theorem}Let $(D,v) \in \mathcal{D}(F)$ be a degree $p^n$ semi-ramified division algebra with $\oD/\oF$ a separable extension which is not strongly degenerate (resp. not degenerate) as in Definition \ref{def1}.  Then for $E/F$ any extension of degree prime to $p$, $D_E$ is a semi-ramified division algebra with $\overline{D_E}/\oE$ a separable extension and $D_E$ is not strongly degenerate (resp. not degenerate).
\label{theorem19}
\end{theorem}

\proof  Let $\theta_D:\Gamma_D/\Gamma_F \to \gal(\overline{D}/\overline{F})$ be the isomorphism given in \eqnref{graded6}.  Let  $\{\pi_i\}_{i=1}^r$ be a set of elements in $D$ such that $\{\theta_D(v(\pi_i)+\Gamma_F)=\sigma_i\}_{i=1}^r$ is a basis for $G\cong \gal(\overline{D}/\overline{F})$.  The hypothesis that $D$ is not strongly degenerate (resp. not degenerate) implies that the matrix $\overline{u}$ is not strongly degenerate (resp. not degenerate) in $\oD/\oF$
.  By Lemma \ref{lemma34}, $D_E$ is semi-ramified and $\overline{D_E}/\oE$ is a separable extension.  Moreover, $\overline{D_E} \cong \overline{D}\otimes_{\oF} \overline{E}$ and is the field join of $\oD$ and $\oE$ and the homomorphism 
\[\theta_{D_E}:\Gamma_{D_E}/\Gamma_E \to \gal(\overline{D_E}/\overline{E})\]
is an isomorphism.  Consider the set of elements $\theta_{D_E}(v(\pi_i \otimes 1)+\Gamma_E)\in \gal(\overline{D_E}/\overline{E})$.  Since $\theta_{D_E}(v(\pi_i \otimes 1)+\Gamma_E)(\overline{d})=\sigma_i(\overline{d})$ for all $\overline{d} \in \oD\subseteq \overline{D_E}$, $\theta_{D_E}(v(\pi_i \otimes 1)+\Gamma_E)=\widetilde{\sigma_i}$, where $\sigma_i \mapsto \widetilde{\sigma_i}$ under the isomorphism $\gal(\overline{D}/\overline{F}) \cong \gal(\overline{D}\,\overline{E}/\overline{E})$.  In particular, $\{\theta_{D_E}(v(\pi_i \otimes 1)+\Gamma_E)\}_{i=1}^r$ is a basis of $G \cong \gal(\overline{D_E}/\overline{E})$.  Since
\[(\pi_i\otimes1)(\pi_j\otimes1)(\pi_i\otimes1)^{-1}(\pi_j\otimes1)^{-1}=u_{ij}\otimes1\in D_E\]
and $\overline{u_{ij}\otimes 1}$ is the image of $\overline{u_{ij}}$ in the inclusion $\oD \hookrightarrow \overline{D_E}$, it is enough to show that the matrix $\overline{u}=(\overline{u}_{ij})$ is not strongly degenerate (resp. not degenerate) in $\overline{D_E}$.  This follows directly from Theorem \ref{theorem8} since $\ou$ is not strongly degenerate (resp. not degenerate) in $\oD$ and $\overline{D_E}$ is obtained from $\oD$ by extending scalars to $\oE$, a prime to $p$ extension of $\oF$. \endproof

\section{The generic abelian crossed product}\label{sec3}
\subsection{Prime to $p$ extensions of generic abelian crossed products}
Let $K/F$, be a noncyclic abelian Galois extension with Galois group $G=\G$ and let $\Delta = \left( K/F,w_{\sigma},u,b\right)$ be an abelian crossed product.  Given $x_1, \ldots, x_r$, independent indeterminates over $K$, consider $K(x_1, \ldots, x_r)/F(x_1, \ldots, x_r)$, the Galois extension induced from $K/F$.  One easily checks that the set of elements $\{u_{ij},b_ix_i\}_{i,j=1}^r$ satisfy relations (\ref{equation1})-(\ref{equation4}) in the Galois extension $K(x_1, \ldots, x_r)/F(x_1, \ldots, x_r)$ and therefore, by \cite{AS} Theorem 1.3, we can make the following definition.
\begin{definition}[\cite{AS} Theorem 2.3]  Let $bx = \{b_i x_i\}_{i=1}^r$.  The crossed product 
\[ \mathcal{A}_{\Delta} = \left(K(x_1, \ldots, x_r)/F(x_1, \ldots, x_r), z_{\sigma},u, bx \right)\]
is the {\bf generic abelian crossed product} defined by $\Delta$.
\label{defn4}
\end{definition}
Throughout this chapter we will set $K'=K(x_1,\ldots,x_r)$ and $F'=F(x_1,\ldots,x_r)$.  In \cite{AS}, Theorem 2.3, $\mathcal{A}_{\Delta}$ is shown to be a division algebra isomorphic to the central localization of an iterated twisted polynomial ring. In particular, $\mathcal{A}_{\Delta} \cong \kssu$ where  $\kssu=K(s_1,...,s_r;\sigma_1,...,\sigma_r;u)$
is the central localization of the twisted polynomial ring 
\[K[s;\sigma;u]=K[s_1,...,s_r;\sigma_1,...\sigma_r;u],\]  and the ring $K[s;\sigma;u]$ is the polynomial ring $K[s_1,\ldots,s_r]$ as a set, with multiplication twisted by $s_i\alpha = \sigma_i(\alpha)s_i$ for $\alpha \in K$ and $s_is_j=u_{ij}s_js_i$.  The isomorphism $\mathcal{A}_{\Delta} \cong \kssu$ is the identity on $K$ and maps $z_i \mapsto s_i$ (the details are found in \cite{AS}).  The center of $K[s;\sigma;u]$ is $F[X_1,...,X_r]$ where $b_i^{-1}s_i^{n_i}=X_i$ and one can check that $x_i \mapsto X_i$ in the isomorphism $\mathcal{A}_{\Delta} \cong \kssu$.  For any element $t \in K[s;\sigma;u]$, let $t^v$ be the unique monomial in $t$ with smallest degree where $\N^r$ is ordered with respect to right to left lexicographical ordering.  This function satisfies $(t^v)^q=(t^q)^v$ for all nonnegative integers $q$.

The next lemma is key in knowing precisely when generic abelian crossed product $p$-algebras have $p$-power central elements.
\begin{lemma}Let $\mathcal{A}_{\Delta}=(K'/F',z_{\sigma},u,bx)$ be the generic abelian crossed product associated to  
\[\Delta=(K/F,w_{\sigma},u,b),\] 
 with $\chara(F) = p>0$ and $G=\G$ a noncyclic $p$-group.  $\mathcal{A}_{\Delta}$ has a $p$-power central element if and only if $u$ is strongly degenerate in $K^*$.\label{lemma21}
\end{lemma} 
\proof  Let $y \in \kssu$ be a nontrivial $p$-power central element.  Then $y=t/t_0$ for some $t \in K[s;\sigma;u]$ and $t_0 \in Z(K[s;\sigma;u])=F(X_1,\ldots,X_r)$.  The element $y$ is $p$-power central if and only if $t$ is $p$-power central.  Let $t_c$ be the sum of monomials in $t$ which are central.  Since $t_c$ is central, $(t-t_c)^p = t^p-t_c^p\in F[x_1,...,x_r]$, i.e., $t-t_c$ is $p$-power central.  Therefore, without loss of generality we may assume that no monomial in $t$ is central.  Moreover, since $(t^v)^p= (t^p)^v$, the monomial $(t^v)$ is $p$-power central.  Let $t^v = l's^{\barm'}$, a monomial in $K[s;\sigma;u]$ with coefficient $l' \in K$.  In $s^{\barm'} = s_1^{m_1'}\ldots s_r^{m_r'}$, it is possible that $m_i'\geq n_i$.  Substituting in $b_iX_i$ for each $s_i^{n_i}$, we see that there exists an $l \in L^*$, $\barm \in \N^r$ with $m_i \in \{0,\ldots,n_i-1\}$ and $\barw \in \N^r$ such that $t^v =ls^{\barm}X^{\barw}$.  Since $X^{\barw}$ is central, the monomial $ls^{\barm}$ is $p$-power central.  Moreover $K/F$ is Galois and hence a separable extension and therefore $\barm \ne \overline{0}$.  Since $(ls^{\barm})^p$ is central it commutes with $K$ and therefore $(\sigma^{\barm})^p=e$, the identity element of $G$.  

We have shown that $\kssu$ has a $p$-power central element if and only if it has a $p$-power central monomial of the form $ls^{\barm}$ with $\sigma^{\barm}$ having order $p$ in $G$.  Under the isomorphism $\kssu\cong \mathcal{A}_{\Delta} $, $ls^{\barm} \mapsto lz^{\barm}$ with $l \in K^*$.  We now show that $lz^{\barm}$ is $p$-power central in $\mathcal{A}_{\Delta}$ if and only if the homogeneous element $lw^{\barm}$ is $p$-power central in $\Delta$.  Calculating $(lz^{\barm})^p$ we get
\begin{eqnarray}
(lz^{\barm})^p&=&l''z_1^{pm_1}\ldots z_r^{pm_r}\nonumber \\
&=&l''(b_1x_1)^{q_1}\ldots(b_rx_r)^{q_r}\nonumber\\
&=& l''b^{\overline{q}}x^{\overline{q}}\in F\cdot x^{\overline{q}}\subset F(x_1,\ldots, x_r)\label{eqn78}
\end{eqnarray}
where $q_i=\frac{m_i\,p}{n_i}$ and $l'' \in K$.  And, 
\begin{eqnarray} 
(lw^{\barm})^p&=& l''w_1^{m_1p}\ldots w_r^{m_rp}\nonumber\\
&=& l''b^{\overline{q}}.\label{eqn79}
\end{eqnarray}
By \eqnref{eqn78} $l''b \in F$ and therefore by \eqnref{eqn79} $(lw^{\barm})^p = l''b^{\overline{q}} \in F$.  By Lemma \ref{lemma22} this happens if and only if $u$ is strongly degenerate in $K$.  \endproof
Before we state the main result about generic abelian crossed products, we need some preliminary lemmas about the power series version of the generic abelian crossed product.  Let $F''$ be the iterated Laurent series field in the $r$ variables $x_1, \ldots, x_r$.  That is, for $1\leq i \leq r$, set $F_0=F$ and $F_i=F_{i-1}((x_i))$.  Then $F''=F_r$.  We will also use the notation $F''=F((x_1,\ldots,x_r))$.  Similarly set $K''=K((x_1,\ldots, x_r))$.  
\begin{lemma} Let $p$ be a fixed prime and let $E/F'$ be an extension of fields of degree prime to $p$.  Then there exists a composite extension $EF''$ of $E$ and $F''$ over $F'$ which has degree prime to $p$ over $F''$.\label{lemma42}
\end{lemma}
\proof Let $v:F' \to \Z^r$ be the standard valuation with $v(x_i)=\overline{e}_i$.  The completion of $F'$ with respect to this valuation is $F''$.  Assume $E$ is a separable extension over $F'$.  Then $E\otimes_{F'}F''$ is canonically isomorphic to $\prod_{w|v}E_w$ where, for every valuation $w$ on $E$ lying over $v$ on $F$, $E_w$ is the completion of $E$ with respect to $w$.  Since the dimension of $E\otimes_{F'}F''$ over $F''$ is prime to $p$ at least one $E_w$ will have degree prime to $p$ over $F''$.  Let $E_{w_0}$ be one of the $E_w$'s with degree over $F''$ prime to $p$.  Then $E_{w_0}$ is our desired composite extension. 

Now let $E$ be an arbitrary extension of $F'$ of degree prime to $p$.  If $\chara(F')=0$ or $\chara(F')=p$, then $E$ is a separable extension and we are done.  So assume $\chara(F')=q$ with $q\ne p$ and $q \ne 0$.  Let $E_s$ be the separable closure of $F'$ in $E$ and let $E_sF''$ be a composite extension of $E_s$ and $F''$ of degree prime to $p$ over $F''$.  Set $E=E_s(\alpha_1,\ldots,\alpha_n)$.  Any composite extension of $E$ and $E_sF''$ is given by an embedding of the $\alpha_i$ into an algebraic closure of $F''$.  Any such embedding will take purely inseparable elements to purely inseparable elements, implying that a composite, $EF''$ of $E$ and $E_sF''$ will be a purely inseparable extension.  Therefore $n_i$, the degree of $EF''$ over $E_sF''$, will be a power of $q$ which is prime to $p$.  Therefore $EF''$ is a composite extension of $E$ and $F''$ over $F'$ and has degree prime to $p$ over $F''$.
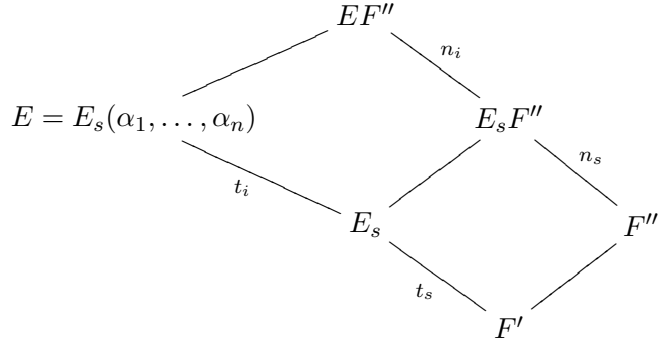
\[\xymatrix{
                        &       EF'' \ar@{-}[dl]        \ar@{-}[dr]     ^{n_i}& &\\
E=E_s(\alpha_1,\ldots,\alpha_n)\ar@{-}[dr]_{t_i}      &                                       &E_sF''\ar@{-}[dl]\ar@{-}[dr]^{n_s}&\\
                        &E_s \ar@{-}[dr]_{t_s} & &F'' \ar@{-}[dl]\\
                        &                       &F'&
}\]    
\endproof
As in the rational field case of Definition \ref{defn4}, given an abelian crossed product $\Delta = \left( K/F,w_{\sigma},u,b\right)$ the elements $\{u_{ij},b_ix_i\}_{i,j=1}^r$ satisfy relations (\ref{equation1})-(\ref{equation4}) in the abelian Galois extension $K((x_1, \ldots, x_r))/F((x_1, \ldots, x_r))$.  Therefore, by \cite{AS} Theorem 1.3, the following definition can be made.
\begin{definition} The abelian crossed product 
\[ \A_{\Delta} = \left(K((x_1, \ldots, x_r))/F((x_1, \ldots, x_r)), z_{\sigma},u, bx \right)\]
is the {\bf power series generic abelian crossed product} defined by $\Delta$.  
\label{defn}
\end{definition}
\begin{remark}This algebra is also defined in \cite{Tignol}.\end{remark}
\begin{lemma}
Let $\Delta=(K/F,z_{\sigma},u,b)$ be an abelian crossed product with $G=\G$, $r \geq 2$, and defining matrix $u$ not strongly degenerate (resp. non-degenerate).  Let $\A_{\Delta}$ be the power series generic abelian crossed product defined by $\Delta$.  Then $\A_{\Delta}$ is a not strongly degenerate (resp. non-degenerate) division algebra as in Definition \ref{def1}.\label{lemma40}
\end{lemma}
\proof  By Definition \ref{def1} we need to show that $\A_{\Delta}$ is a semi-ramified division algebra with separable residue fields and that there exists a choice of elements which push forward to a basis of $G$ under $\theta_{\A_{\Delta}}$ and define a not strongly degenerate (resp. non-degenerate) matrix in $\overline{\A_{\Delta}}/\overline{F''}$.  Let $v:K'' \to \Z^r$ be the standard Henselian valuation with $v(x_i)=\overline{e}_i$.  With respect to this valuation, the residue fields of $K''$ and $F''$ are $\overline{K''}\cong K$ and $\overline{F''}\cong F$.  Since $[K'':F'']=[\overline{K''}:\overline{F''}]$, and $K/F$ is a separable extension, the subfield $K''$ is an inertial maximal subfield of $\A_{\Delta}$.  Therefore, by \cite{JW} Lemma 5.1, $Z(\overline{\A_{\Delta}})$ is separable over $F$, $\theta_{\A_{\Delta}}$ is an isomorphism, and $\A_{\Delta}$ is defectless over $F''$.  Since $v(z_i) = \frac{1}{n_i}v(x_i)$, the index of $\Gamma_{F''}$ in $\Gamma_{\A_{\Delta}}$ is $\geq |G|$.  Moreover, $[\overline{\A_{\Delta}}:\overline{F''}]\geq [\overline{K''}:\overline{F''}]=|G|$.  Since 
\[|G|^2 = [\A_{\Delta}:F''] = [\overline{\A_{\Delta}}:\overline{F''}] [\Gamma_{\A_{\Delta}}: \Gamma_{F''}],\] 
these two inequalities are equalities, that is $\overline{\A_{\Delta}} \cong K$, $[\overline{\A_{\Delta}}:\overline{F''}]=|G|$ and $[\Gamma_{\A_{\Delta}}:\Gamma_{F''}]=|G|$.  It is now clear that $\A_{\Delta}$ is a semi-ramified division algebra with $\overline{\A_{\Delta}}/\overline{F''}$ a separable extension.  Since inner automorphism by $z_i$ in $\A_{\Delta}$ induces $\sigma_i$ on $K''$ we have $\theta_{\A_{\Delta}}(v(z_i)+\Gamma_{F''}) = \sigma_i$.  By hypothesis the matrix defined by  $u_{ij} = z_iz_jz_i^{-1}z_j^{-1}$ is not strongly degenerate (resp. non-degenerate) in $\overline{\A_{\Delta}}\cong K$, therefore, $\A_{\Delta}$ is a not strongly degenerate (resp. non-degenerate) division algebra.   \endproof

\begin{theorem} Let $\mathcal{A}_{\Delta}$ be the generic abelian crossed product associated to  
\[\Delta=(K/F,w_{\sigma},u,b),\] 
 with $\chara(F) = p>0$ and $G$ a noncyclic abelian $p$-group.  The matrix $u$ is strongly degenerate if and only if there exists a prime to $p$ extension $E/F'$ such that $\mathcal{A}_{\Delta} \otimes_{F'} E$ contains a nontrivial $p$ power central element.  
\label{graded4}
\end{theorem}
\proof  Assume $u$ is strongly degenerate in $K$.  Then, by Lemma \ref{lemma21} $\cA_{\Delta}$ has a nontrivial $p$-power central element and so we are done with one direction.  

On the other hand assume there exist $E/F'$ a prime to $p$ extension such that $\cA_{\Delta}$ has a $p$-power central element.  Assume by way of contradiction that $u$ is not strongly degenerate.  Let $\A_{\Delta}$ be the power series generic abelian crossed product associated to $\Delta$.  By Lemma \ref{lemma40} $\A_{\Delta}$ is a not strongly degenerate division algebra.  Let $E''=EF''$ be a composite extension of $E$ and $F''$ with degree prime to $p$ over $F''$.  Such an extension exists by Lemma \ref{lemma42}.  By Lemma \ref{lemma34}, $\A_{\Delta} \otimes_{F''} E''$ is a semi-ramified inertially split division algebra and, by Theorem \ref{theorem19}, $\A_{\Delta}\otimes E''$ is not strongly degenerate.  By Corollary \ref{corollary10} $\A_{\Delta}\otimes E$ has no $p$-power central elements.  This is a contradiction since we assumed $\mathcal{A}_{\Delta}\otimes_{F'}E$ has a non-trivial $p$-power central element which implies $\A_{\Delta}\otimes E''$ has one as well.  Therefore $u$ is strongly degenerate.\endproof

\begin{corollary}Let $\Delta$, $\cA_{\Delta}$ be as in Theorem \ref{graded4}.  If $u$ is not strongly degenerate then for all prime to $p$ extensions $E/F'$, $\cA_{\Delta}\otimes_{F'}E$ has no nontrivial $p$-power central elements.  In particular, $\mathcal{A}_{\Delta}$ does not become cyclic after tensoring by any prime to $p$ extension. 
\label{cor89}
\end{corollary}
\proof This is simply stating one contrapositive of Theorem \ref{graded4}.\endproof
\begin{remark}  Since a matrix which is not degenerate is not strongly degenerate, Corollary \ref{cor89} shows that the non-cyclic $p$-algebras constructed by Amitsur and Saltman in \cite{AS}, Theorem 3.2, remain non-cyclic after tensoring by any prime to $p$ extension.
\end{remark}
\subsection{Indecomposable generic abelian crossed products}
In this section we prove that the indecomposable division algebras constructed by Saltman in \cite{LN}, Theorem 7.17, remain indecomposable after any prime to $p$ extension.  We start by proving a preliminary lemma which is independent of the characteristic of the ground field.
\begin{lemma}
Let $p$ be a fixed prime and let $\mathcal{A}_{\Delta}$ be a generic abelian crossed product.  If $\mathcal{A}_{\Delta}$ is decomposable after a prime to $p$ extension then the power series generic abelian crossed product $\A_{\Delta}$ defined by the same crossed product $\Delta$ is decomposable after a prime to $p$ extension.  \label{lemma43}
\end{lemma}

\proof  Let $F'=F(x_1,\ldots,x_r)$ be the center of $\mathcal{A}_{\Delta}$ and $F'' = F((x_1,\ldots,x_r))$ the center of $\A_{\Delta}$.  Assume there exists an extension $E/F'$ of degree prime to $p$ such that $\mathcal{A}_{\Delta} \otimes_{F'}E$ is decomposable.  Set $\mathcal{A}_{\Delta} \otimes_{F'}E \cong D_1 \otimes_{E}D_2$ with $\deg(D_i)>1$.  By Lemma \ref{lemma42} there exists a composite extension, $EF''$, which has degree over $F''$ prime to $p$.  Then, 
\begin{equation}\label{eqn68}
\begin{split}
(\mathcal{A}_{\Delta} \otimes_{F'}E)\otimes_E EF''&=\mathcal{A}_{\Delta} \otimes_{F'}(E\otimes_E EF'')\\
&=\mathcal{A}_{\Delta} \otimes_{F'} (F''\otimes_{F''}EF'')\\
&=(\mathcal{A}_{\Delta} \otimes_{F'} F'') \otimes_{F''}EF''\\
&=\A_{\Delta}\otimes_{F''} EF''.
\end{split}
\end{equation}
Furthermore,
\begin{equation}\label{eqn69}
(\mathcal{A}_{\Delta}\otimes_{F'}E) \otimes_{E}EF'' \cong (D_1 \otimes_{F'}EF'') \otimes_{EF''}(D_2 \otimes_{F'}EF'').
\end{equation}
Combining \eqnref{eqn68} and \eqnref{eqn69} we see that $\A_{\Delta}$ is decomposable after a prime to $p$ extension.  \endproof

In the case $\gal(K/F)=G=C_p\times C_p$, the elementary abelian group of order $p^2$, the condition $u$ is not degenerate takes on a form which is comparatively easy to describe.  Let $\ig K^*$ be the multiplicative subgroup of $K^*$ generated by elements of the form $\sigma(x)\inv{x}$ for $x \in K^*$ and $\sigma \in G$.
\begin{lemma}Let $K/F$ be a noncyclic elementary abelian extension with $\gal(K/F)=\langle\sigma_1\rangle\times\langle\sigma_2\rangle$, $|\sigma_i|=p$.  Let $\Delta=(K/F,z_{\sigma},u,b)$ be an abelian crossed product.  The matrix $u$ is degenerate if and only if $u_{12} \in \ig K^*$.\label{rank2degenerate}
\end{lemma}
\proof  Assume $u$ is degenerate, that is, there exist $\sigma^{\barm},\sigma^{\barn} \in G$ such that $\langle \sigma^{\barm},\sigma^{\barn}\rangle$ is not cyclic and
\[u_{\barm,\barn} = \frac{\sigma^{\barm}(a)}{a}\frac{\sigma^{\barn}(b)}{b}\]
for some $a,b, \in K^*$.  Since $\abs{G}=p^2$ this implies that $\langle \sigma^{\barm},\sigma^{\barn}\rangle=G$.  Choose non-negative integers $s_1,t_1,s_2,t_2 \in \N$ such that $(\sigma^{\barm})^{s_i}(\sigma^{\barn})^{t_i}=\sigma_i$.  The degeneracy condition implies $bz^{\barm}$ and $a^{-1}z^{\barn}$ commute.  Therefore, the elements $(bz^{\barm})^{s_i}(a^{-1}z^{\barn})^{t_i}$  for $i=1,2$ also commute and conjugation by these elements induce the action of $\sigma_i$ on $K$.  By \cite{LN} Proposition 7.13, this implies $u_{12} \in \ig K^*$.  Conversely, assume $u_{12} \in \ig K^*$.  Since $\ig$, the augmentation ideal of the group ring $\zg$, is generated by $\sigma_1-1$ and $\sigma_2-1$ we in fact have $u_{12}=\sigma_1(a)a^{-1}\sigma_2(b)b^{-1}$ for some $a,b \in K^*$.  In particular, $u$ is degenerate.\endproof

\begin{theorem}Let $G=\langle \sigma_1\rangle \times \langle \sigma_2 \rangle$ with $|\sigma_i|=p$, $p\ne 2$.  Let $\Delta = (K/F,z_{\sigma},u,b)$ be an abelian crossed product with exponent $p$ and $u_{12} \notin \ig K^*$.  Then $\mathcal{A}_{\Delta}$ is an indecomposable division algebra of exponent $p$ and index $p^2$.  Moreover, $\cA_{\Delta}$ remains indecomposable after any prime to $p$ extension.
\label{theorem23}
\end{theorem}
\begin{remark} In \cite{LN}, Theorem 7.17, Saltman proves that the generic abelian crossed products in Theorem \ref{theorem23} are indecomposable with the added assumption that $F$ contains a primitive $p$-th root of unity.
\end{remark}
\begin{remark}  The condition $p \ne 2$ in Theorem \ref{theorem23} is not an unnecessary restriction.  This follows because any division algebra of exponent 2 and index 4 is biquaternion by \cite{Albert}, Theorem 11.2.\end{remark}
\proof  $\mathcal{A}_{\Delta}$ is index $p^2$ and exponent $p$ by \cite{LN}, Theorem 7.17.  Assume $\chara(F)=p$ and $\cA_{\Delta}\otimes_{F'}E\cong D_1 \otimes_E D_2$ is a nontrivial decomposition of $\cA_{\Delta}\otimes_{F'}E$ where $E$ is a prime to $p$ extension of $F'$.  Since the $D_i$ have prime index $p$, after a further prime to $p$ extension, at least one of the $D_i$ become cyclic.  In particular, after a prime to $p$ extension, $\cA_{\Delta}$ has a $p$-power central element.  This is a contradiction to Theorem \ref{graded4}.

Assume $\chara(F) \ne p$.  It is enough to show, by Lemma \ref{lemma43}, that $\A_{\Delta}$, the power series generic abelian crossed product, does not become decomposable after a prime to $p$ extension.  Assume that there is a prime to $p$ extension $E/F''$ and a non-trivial decomposition $\A_{\Delta} \otimes_{F''} E \cong D_1 \otimes_E D_2$ with $\deg(D_i)=p$.  By Remark \ref{rank2degenerate} the condition $u \notin \ig K^*$ is equivalent to the condition that the matrix $u$ is non-degenerate.  Therefore, by Lemma \ref{lemma40}, $\A_{\Delta}$ is a non-degenerate division algebra and in particular it is not strongly degenerate. Let $v$ denote the extension of the standard Henselian valuation on $F''$ to $\A_{\Delta}\otimes_{F''}E$.  As in Theorem \ref{graded4} $(\A_{\Delta}\otimes_{F''}E,v)$ is a semi-ramified division algebra with separable residue field which is not strongly degenerate.  We can therefore apply Theorem \ref{graded1} to study the value of non-central elements in $D_1$ and $D_2$.

Let $\alpha \in D_1-E$.  Then, since $D_1$ has degree $p$ over $E$, $\alpha$ satisfies a polynomial relation
\begin{equation}\alpha^p+s_1\alpha^{p-1}+\ldots +s_{p-1}\alpha = s_p \label{graded5}\end{equation}
with each $s_i \in E$.  We will show $v(\alpha) \in \Gamma_E$.  Since $[E(\alpha):E]=p$ we have either $\abs{\Gamma_{E(\alpha)}:\Gamma_E} = p$ or $\abs{\Gamma_{E(\alpha)}:\Gamma_E} =1$.   In either case $v(\alpha^p) \in \Gamma_E$.  Set $D=\A_{\Delta}\otimes_{F''}E$ and let $GD$ be the graded division algebra associated to $(D,v)$.  Since $D$ is not strongly degenerate, by Theorem \ref{graded1} all $p$-power central homogeneous elements in $GD$ are in $D_{\gamma}$ for $\gamma \in \Gamma_E$.  That is, if we show $g(\alpha^p) \in GE=Z(GD)$ then by \ref{graded1} $v(\alpha) \in \Gamma_E$.   From equation \eqnref{graded5} and the triangle inequality, we see that to show that $g(\alpha^p)\in GE$ it is enough to show that $v(s_i\alpha^{p-i})>v(\alpha^p)$ for each $1\leq i \leq p-1$.  If $v(\alpha) \notin \Gamma_E$ then $\abs{\Gamma_{E(\alpha)}:\Gamma_E}=p$ and the elements $\{v(s_i\alpha^{p-i})\}_{i=1}^{p-1}$ are all distinct in $\Gamma_{E(\alpha)}-\Gamma_E$.  Hence 
\[v(s_p) = \textrm{min}\{v(\alpha^p),v(s_i\alpha^{p-i})\}_{i=1}^{p-1}=v(\alpha^p)\in \Gamma_E. \]
Therefore, $v(\alpha^p)<v(s_i\alpha^{p-i})$ for each $i=1,\ldots,p-1$ and hence $g(\alpha^p)=g(s_p) \in GE$.  Therefore by Theorem \ref{graded1} $v(\alpha) \in \Gamma_E$.  Since $\alpha$ was an arbitrary element of $D_1$, every element of $D_1$ has value in $\Gamma_E$.  
And, since this argument was not particular to $D_1$, $D_2$ also has every element having valuation in $\Gamma_E$.  Since $\overline{F''}\cong F$ has characteristic not $p$, both $D_1$ and $D_2$ are defectless.  Therefore, for $i=1,2$, $p^2=[D_i:E]=[\overline{D_i}:\overline{E}]\abs{\Gamma_{D_i}:\Gamma_E}=[\overline{D_i}:\overline{E}]$.  Since $D$ is semi-ramified, $[\overline{D}:\overline{E}]=p^2$.  Therefore, as $\overline{D_i}\subseteq\overline{D}$ and they both have dimension $p^2$ over $\overline{E}$, $\overline{D_i} \cong \overline{D}$.  This implies that $\oD_i$ is a separable field extension of $E$.  That is, the $D$-subdivision algebras $D_i$ are inertial over $E$.  But, by \cite{JW}, Lemma 2.2, if $D_i$ is inertial over $E$ then $\overline{E}=\overline{Z(D_i)}=Z(\overline{D_i})=\overline{D_i}$.  This contradicts that $[\overline{D_i}:\overline{E}]=p^2$.\endproof

\begin{remark}The existence of abelian crossed products $\Delta$ satisfying the conditions of Theorem \ref{theorem23} is addressed in \cite{LN}.  In \cite{LN}, Theorem 12.14, an abelian crossed product $\Delta$ is constructed by taking a generic $C_p \times C_p$-crossed product with exponent equal index equal $p^2$ and tensoring by the Severi-Brauer splitting field of the $p$-th tensor power.  It is important to mention that this particular $\Delta$ is an example of a ``generic algebra'' as defined in \cite{Karpenko2}.  Moreover, by \cite{Karpenko2}, Corollary 5.4, $\Delta$ is indecomposable and remains indecomposable after any prime to $p$ extension.  Karpenko's results were obtained by studying torsion in $CH^2$ of Severi-Brauer varieties.  \cite{LN} Theorem 7.17 and Theorem \ref{theorem23} from above show that in a specific example the generic abelian crossed product associated to these algebras are indecomposable and remain so after any prime to $p$ extension.  Both proofs make no mention of the geometry of Severi-Brauer varieties.  It is not always the case that an indecomposable generic abelian crossed product is constructed from an indecomposable crossed product.  The author plans to prepare a paper discussing an example of this soon.
\end{remark} 
\bibliographystyle{alpha} 
\bibliography{ref}
\end{document}